\documentclass[11pt]{article}

\catcode`\@=11

\thicklines
\newskip\Einheit \Einheit=.6cm
\newcount\xcoord \newcount\ycoord
\newdimen\xdim \newdimen\ydim \newdimen\PfadD@cke \newdimen\Pfadd@cke
\def\PfadDicke#1{\PfadD@cke#1 \divide\PfadD@cke by2 
\Pfadd@cke\PfadD@cke \multiply\PfadD@cke by2}
\long\def\LOOP#1\REPEAT{\def\BODY{#1}\ITERATE}
\def\ITERATE{\BODY \let\next\ITERATE \else\let\next\relax\fi \next}
\let\REPEAT=\fi
\def\Punkt{\hbox{\raise-2pt\hbox to0pt{\hss\scriptsize$\bullet$\hss}}}

\def\DuennPunkt(#1,#2){\unskip
  \raise#2 \Einheit\hbox to0pt{\hskip#1 \Einheit
          \raise-1.5pt\hbox to0pt{\hss\tiny$\bullet$\hss}\hss}}
		  
\def\NormalPunkt(#1,#2){\unskip
  \raise#2 \Einheit\hbox to0pt{\hskip#1 \Einheit
          \raise-3pt\hbox to0pt{\hss\large$\bullet$\hss}\hss}}
\def\DickPunkt(#1,#2){\unskip
  \raise#2 \Einheit\hbox to0pt{\hskip#1 \Einheit
          \raise-4pt\hbox to0pt{\hss\Large$\bullet$\hss}\hss}}
\def\Kreis(#1,#2){\unskip
  \raise#2 \Einheit\hbox to0pt{\hskip#1 \Einheit
          \raise-4pt\hbox to0pt{\hss\Large$\circ$\hss}\hss}}
\def\Diagonale(#1,#2)#3{\unskip\leavevmode
  \xcoord#1\relax \ycoord#2\relax
      \raise\ycoord \Einheit\hbox to0pt{\hskip\xcoord \Einheit
         \unitlength\Einheit
         \line(1,1){#3}\hss}}
\def\AntiDiagonale(#1,#2)#3{\unskip\leavevmode
  \xcoord#1\relax \ycoord#2\relax \advance\xcoord by -0.05\relax
      \raise\ycoord \Einheit\hbox to0pt{\hskip\xcoord \Einheit
         \unitlength\Einheit
         \line(1,-1){#3}\hss}}
\def\Pfad(#1,#2),#3\endPfad{\unskip\leavevmode
  \xcoord#1 \ycoord#2 \thicklines\ZeichnePfad#3\endPfad\thinlines}
\def\ZeichnePfad#1{\ifx#1\endPfad\let\next\relax
  \else\let\next\ZeichnePfad
    \ifnum#1=1
      \raise\ycoord \Einheit\hbox to0pt{\hskip\xcoord \Einheit
         \vrule height\Pfadd@cke width1 \Einheit depth\Pfadd@cke\hss}%
      \advance\xcoord by 1
     \else\ifnum#1=2
      \raise\ycoord \Einheit\hbox to0pt{\hskip\xcoord \Einheit
         \unitlength\Einheit
         \line(0,1){1}\hss}
      \advance\xcoord by 0
      \advance\ycoord by 1
 \else\ifnum#1=3
      \raise\ycoord \Einheit\hbox to0pt{\hskip\xcoord \Einheit
         \unitlength\Einheit
         \line(1,1){1}\hss}
      \advance\xcoord by 1
      \advance\ycoord by 1
    \else\ifnum#1=4
      \raise\ycoord \Einheit\hbox to0pt{\hskip\xcoord \Einheit
         \unitlength\Einheit
         \line(1,-1){1}\hss}
      \advance\xcoord by 1
      \advance\ycoord by -1
   \else\ifnum#1=5
      \raise\ycoord \Einheit\hbox to0pt{\hskip\xcoord \Einheit
         \unitlength\Einheit
         \line(2,1){2}\hss}
      \advance\xcoord by 2
      \advance\ycoord by 1
	  \else\ifnum#1=6
      \raise\ycoord \Einheit\hbox to0pt{\hskip\xcoord \Einheit
         \unitlength\Einheit
         \line(2,-1){2}\hss}
      \advance\xcoord by 2
      \advance\ycoord by -1
	  \else\ifnum#1=7
      \raise\ycoord \Einheit\hbox to0pt{\hskip\xcoord \Einheit
         \unitlength\Einheit
         \line(3,1){3}\hss}
      \advance\xcoord by 3
      \advance\ycoord by 1
	  \else\ifnum#1=8
      \raise\ycoord \Einheit\hbox to0pt{\hskip\xcoord \Einheit
         \unitlength\Einheit
         \line(3,-1){3}\hss}
      \advance\xcoord by 3
      \advance\ycoord by -1
    \fi\fi\fi\fi\fi\fi\fi\fi
  \fi\next}
\def\hSSchritt{\leavevmode\raise-.4pt\hbox 
to0pt{\hss.\hss}\hskip.2\Einheit
  \raise-.4pt\hbox to0pt{\hss.\hss}\hskip.2\Einheit
  \raise-.4pt\hbox to0pt{\hss.\hss}\hskip.2\Einheit
  \raise-.4pt\hbox to0pt{\hss.\hss}\hskip.2\Einheit
  \raise-.4pt\hbox to0pt{\hss.\hss}\hskip.2\Einheit}
\def\vSSchritt{\vbox{\baselineskip.2\Einheit\lineskiplimit0pt
\hbox{.}\hbox{.}\hbox{.}\hbox{.}\hbox{.}}}
\def\DSSchritt{\leavevmode\raise-.4pt\hbox to0pt{%
  \hbox to0pt{\hss.\hss}\hskip.2\Einheit
  \raise.2\Einheit\hbox to0pt{\hss.\hss}\hskip.2\Einheit
  \raise.4\Einheit\hbox to0pt{\hss.\hss}\hskip.2\Einheit
  \raise.6\Einheit\hbox to0pt{\hss.\hss}\hskip.2\Einheit
  \raise.8\Einheit\hbox to0pt{\hss.\hss}\hss}}
\def\dSSchritt{\leavevmode\raise-.4pt\hbox to0pt{%
  \hbox to0pt{\hss.\hss}\hskip.2\Einheit
  \raise-.2\Einheit\hbox to0pt{\hss.\hss}\hskip.2\Einheit
  \raise-.4\Einheit\hbox to0pt{\hss.\hss}\hskip.2\Einheit
  \raise-.6\Einheit\hbox to0pt{\hss.\hss}\hskip.2\Einheit
  \raise-.8\Einheit\hbox to0pt{\hss.\hss}\hss}}
\def\SPfad(#1,#2),#3\endSPfad{\unskip\leavevmode
  \xcoord#1 \ycoord#2 \ZeichneSPfad#3\endSPfad}
\def\ZeichneSPfad#1{\ifx#1\endSPfad\let\next\relax
  \else\let\next\ZeichneSPfad
    \ifnum#1=1
      \raise\ycoord \Einheit\hbox to0pt{\hskip\xcoord \Einheit
         \hSSchritt\hss}%
      \advance\xcoord by 1
    \else\ifnum#1=2
      \raise\ycoord \Einheit\hbox to0pt{\hskip\xcoord \Einheit
        \hbox{\hskip-2pt \vSSchritt}\hss}%
      \advance\ycoord by 1
    \else\ifnum#1=3
      \raise\ycoord \Einheit\hbox to0pt{\hskip\xcoord \Einheit
         \DSSchritt\hss}
      \advance\xcoord by 1
      \advance\ycoord by 1
    \else\ifnum#1=4
      \raise\ycoord \Einheit\hbox to0pt{\hskip\xcoord \Einheit
         \dSSchritt\hss}
      \advance\xcoord by 1
      \advance\ycoord by -1
    \fi\fi\fi\fi
  \fi\next}
\def\Koordinatenachsen(#1,#2){\unskip
 \hbox to0pt{\hskip-.5pt\vrule height#2 \Einheit width.5pt depth1 
\Einheit}%
 \hbox to0pt{\hskip-1 \Einheit \xcoord#1 \advance\xcoord by1
    \vrule height0.25pt width\xcoord \Einheit depth0.25pt\hss}}
\def\Koordinatenachsen(#1,#2)(#3,#4){\unskip
 \hbox to0pt{\hskip-.5pt \ycoord-#4 \advance\ycoord by1
    \vrule height#2 \Einheit width.5pt depth\ycoord \Einheit}%
 \hbox to0pt{\hskip-1 \Einheit \hskip#3\Einheit 
    \xcoord#1 \advance\xcoord by1 \advance\xcoord by-#3 
    \vrule height0.25pt width\xcoord \Einheit depth0.25pt\hss}}
\def\Gitter(#1,#2){\unskip \xcoord0 \ycoord0 \leavevmode
  \LOOP\ifnum\ycoord<#2
    \loop\ifnum\xcoord<#1
      \raise\ycoord \Einheit\hbox to0pt{\hskip\xcoord 
\Einheit\Punkt\hss}%
      \advance\xcoord by1
    \repeat
    \xcoord0
    \advance\ycoord by1
  \REPEAT}
\def\Gitter(#1,#2)(#3,#4){\unskip \xcoord#3 \ycoord#4 \leavevmode
  \LOOP\ifnum\ycoord<#2
    \loop\ifnum\xcoord<#1
      \raise\ycoord \Einheit\hbox to0pt{\hskip\xcoord 
\Einheit\Punkt\hss}%
      \advance\xcoord by1
    \repeat
    \xcoord#3
    \advance\ycoord by1
  \REPEAT}
\def\Label#1#2(#3,#4){\unskip \xdim#3 \Einheit \ydim#4 \Einheit
  \def\lo{\advance\xdim by-.5 \Einheit \advance\ydim by.5 \Einheit}%
  \def\llo{\advance\xdim by-.25cm \advance\ydim by.5 \Einheit}%
  \def\loo{\advance\xdim by-.5 \Einheit \advance\ydim by.25cm}%
  \def\o{\advance\ydim by.25cm}%
  \def\ro{\advance\xdim by.5 \Einheit \advance\ydim by.5 \Einheit}%
  \def\rro{\advance\xdim by.25cm \advance\ydim by.5 \Einheit}%
  \def\roo{\advance\xdim by.5 \Einheit \advance\ydim by.25cm}%
  \def\l{\advance\xdim by-.30cm}%
  \def\r{\advance\xdim by.30cm}%
  \def\lu{\advance\xdim by-.5 \Einheit \advance\ydim by-.6 \Einheit}%
  \def\llu{\advance\xdim by-.25cm \advance\ydim by-.6 \Einheit}%
  \def\luu{\advance\xdim by-.5 \Einheit \advance\ydim by-.30cm}%
  \def\u{\advance\ydim by-.30cm}%
  \def\ru{\advance\xdim by.5 \Einheit \advance\ydim by-.6 \Einheit}%
  \def\rru{\advance\xdim by.25cm \advance\ydim by-.6 \Einheit}%
  \def\ruu{\advance\xdim by.5 \Einheit \advance\ydim by-.30cm}%
  #1\raise\ydim\hbox to0pt{\hskip\xdim
     \vbox to0pt{\vss\hbox to0pt{\hss$#2$\hss}\vss}\hss}%
}
\catcode`\@=12

\usepackage{amsmath,amsthm,}
\usepackage{color}
\def\red{\textcolor{red} }
\def\green{\textcolor{green} }
\usepackage{xspace}
\usepackage[colorlinks=true,
linkcolor=green,
filecolor=brown,
citecolor=green]{hyperref}
\def\blue{\textcolor{blue} }

\def\a{\ensuremath{\mathcal A}\xspace}
\def\b{\ensuremath{\mathcal B}\xspace}

\def\eps{\epsilon}

\begin{document}
\newtheorem{lemma}{Lemma}
\newtheorem{theorem}{Theorem}
\newtheorem{prop}{Proposition}
\newtheorem{cor}{Corollary}
%\vspace*{-5mm}

\begin{center}
{\Large
Some Identities for the Catalan and Fine Numbers     \\ 
}
\vspace{10mm}
DAVID CALLAN  \\
Department of Statistics  \\
\vspace*{-2mm}
University of Wisconsin-Madison  \\
\vspace*{-2mm}
Medical Science Center \\
\vspace*{-2mm}
1300 University Ave  \\
\vspace*{-2mm}
Madison, WI \ 53706-1532  \\
{\bf callan@stat.wisc.edu}  \\
\vspace{5mm}
February 15, 2005
\end{center}

\vspace{5mm}

\begin{center}
  \textbf{Abstract}
\end{center}

\vspace*{-5mm}

\noindent We establish combinatorial interpretations 
of several identities for the Catalan 
and Fine numbers and, along the way, we present some new bijections of independent interest.
Briefly, we show 
that $C_{n}=\frac{1}{n+1}\sum_{k}\binom{n+1}{2k+1} \binom{n+k}{k}$  
counts ordered trees on $n$ edges by number of interior vertices adjacent to a 
leaf and $C_{n}=\frac{2}{n+1}\sum_{k}\binom{n+1}{k+2} 
\binom{n-2}{k}$ counts Dyck $n$-paths by number of long interior 
inclines. We also give an analogue for the Fine numbers of Touchard's 
Catalan number identity.

\vspace{8mm}

{\Large \textbf{1 \quad Introduction}  }

The purpose of this paper is to establish combinatorial interpretations 
of several identities for the Catalan numbers and one for the Fine 
numbers.
In the bijective spirit, we eschew generating functions and instead use  
bijections, some old, some new. The common thread is a generalization,
from ordinary lattice paths to ``marked'' paths,
of a well known bijective method for counting Dyck paths. Before reviewing 
terminology and introducing marked paths, let us list the main 
results:
\[
C_{n}=\frac{2}{n+1}\sum_{k=0}^{n-2}\binom{n+1}{k+2} 
\binom{n-2}{k}
\]
counts Dyck $n$-paths by number of long interior inclines,
\[
C_{n}=\frac{1}{n+1}\sum_{k}\binom{n+1}{2k+1} \binom{n+k}{k}
\]
counts ordered trees on $n$ edges by number of interior vertices adjacent to a 
leaf, and
\[
F_{n}=\frac{1}{n+1}\sum_{k\ge 0}\binom{n-2-k}{k}2^{n-2-2k}\binom{n+1}{k+1}
\]
counts Fine $n$-paths by number of long noninitial ascents.

\vspace{7mm}
%\newpage 

{\Large \textbf{2 \quad Terminology and Notation}  }

A \emph{balanced $n$-path} is a sequence of $n$ $U$s and $n\ D$s, 
represented as a path of upsteps $(1,1)$ and 
downsteps $(1,-1)$ from $(0,0)$ to $(2n,0)$.
An ascent in a balanced path is a maximal sequence of (consecutive) upsteps and analogously 
for a descent. An incline is an ascent or descent. An incline is 
short if it consists of just one step, otherwise it is long. 
An ascent consisting of $j$ upsteps contains $j-1$ vertices of the 
path in its interior. We will need the notion of a path with marked 
vertices. Two types of vertex are relevant. An IA vertex is one 
incident with two upsteps, that is, a vertex in the interior of an 
ascent (IA for interior ascent). A DF vertex is one that is not 
incident with a downstep (DF for downstep free). Every IA vertex is 
also DF but the initial and terminal vertices of a balanced path, 
while never IA, may or may not be DF. A $k$-marked IA (resp. DF) balanced path is 
one in which $k$ of the IA (resp. DF) vertices have been marked. 

A \emph{Dyck} path is 
a balanced path that never drops below the $x$-axis (ground level).
The size of a Dyck path, sometimes called its semilength, is the number of 
upsteps; thus a Dyck $n$-path has size $n$. The empty Dyck path is 
denoted $\eps$. A nonempty Dyck path always has an initial ascent and 
a terminal descent; all other inclines are interior. A peak is an 
occurrence of $UD$ and similarly a valley is a $DU$. A $DXD$ is an 
occurrence of $DUD$ or $DDD$. A hill is a peak at level 1. Every $U$ 
in a Dyck path has a matching $D$: the downstep terminating the 
\emph{shortest} Dyck subpath beginning at $U$. We also need the 
notion of an associated $D$ with each $U$: the downstep terminating the 
\emph{longest} Dyck subpath beginning at $U$. 

The marked Dyck path 
illustrated has size 6, 4 peaks, 3 valleys, 2 $UU$s, 2 $DD$s, 2 long 
interior inclines, 3 $DXD$s and 1 hill. It has one marked vertex and 
is both IA and DF.
% 1=step in x-direction, 2=step in y-direction, 3=upward diagonal step, 
% 4=downward diagonal step
\Einheit=0.6cm
\[
\Pfad(-6,0),334343344434\endPfad
\SPfad(-6,0),111111111\endSPfad
\SPfad(4,0),11\endSPfad
\Label\o{\textrm{ {\scriptsize $U_{1}$} }}(-5.8,0.4)
\Label\o{\textrm{ {\scriptsize $U_{2}$} }}(-4.8,1.4)
\Label\o{ \textrm{{\scriptsize associated}}}(3.6,2.6)
\Label\o{ \textrm{{\scriptsize $D$ for $U_{2}$}}}(3.6,2.0)
\Label\o{ \textrm{{\scriptsize $\swarrow$ }}}(3.0,1.4)
\Label\o{ \textrm{{\scriptsize matching}}}(3.5,-1.0)
\Label\o{ \textrm{{\scriptsize $D$ for $U_{1}$}}}(3.5,-1.6)
\Label\o{ \textrm{{\scriptsize $\uparrow$ }}}(3.5,-0.4)
\DuennPunkt(-6,0)
\DuennPunkt(-5,1)
\DuennPunkt(-4,2)
\DuennPunkt(-3,1)
\DuennPunkt(-2,2)
\DuennPunkt(-1,1)
\NormalPunkt(0,2)
\DuennPunkt(1,3)
\DuennPunkt(2,2)
\DuennPunkt(3,1)
\DuennPunkt(4,0)
\DuennPunkt(5,1)
\DuennPunkt(6,0)
\]
\begin{center}\nopagebreak 
a marked Dyck path
\end{center}

%\vspace*{7mm}
\newpage

{\Large \textbf{3 \quad The $\mathbf{\frac{1}{n+1}}$ Factor}  }

It is classic that the parameter $X$ on balanced $n$-paths 
defined by $X=$ ``number of upsteps above ground level'' is uniformly
distributed over $\{0,1,2,\ldots,n\}$ and hence divides the $\binom{2n}{n}$ balanced $n$-paths
into $n+1$ equal-size classes, one of which consists of the Dyck $n$-paths (the one with $X=n$).   
Indeed, for $1\le i \le n$, 
a bijection from balanced $n$-paths with $X=0$ (inverted Dyck paths) to those with $X=i$ is as follows.
Number the upsteps from left to right and top to bottom, starting with the last upstep. Then remove 
the first downstep $D$ encountered directly west of upstep $i$ to obtain two 
subpaths $P$ and $Q$, and reassemble as $Q\,D\,P$.
(See \cite{callan} for a more leisurely account.) 

This bijection does not disturb the lengths of ascents. So we can 
conclude that the number of Dyck $n$-paths with $k$ odd-length ascents 
is $\frac{1}{n+1}$ times the number of balanced $n$-paths with $k$ 
odd-length ascents. We will use this fact in \S 8 below.
More importantly, the bijection can 
equally well be applied to $k$-marked IA balanced $n$-paths: the interior 
vertices of ascents are never disturbed. It again produces 
$n+1$ equal-size classes one of which consists of the $k$-marked IA
Dyck $n$-paths. 
One consequence \cite{polydissect} is that the total number of IA  Dyck $n$-paths (with no restriction
on the number of marks) is the little Schr\"{o}der number $s_{n}$ and, 
since DF Dyck paths also allow a mark on the initial vertex, the number 
of DF  Dyck $n$-paths is the big Schr\"{o}der number $r_{n}=2s_{n}\ 
(n\ge 1)$. In fact, there is a simple bijection from DF Dyck 
paths to a well known manifestation of the big Schr\"{o}der 
numbers---Schr\"{o}der paths. A Schr\"{o}der path is a path of upsteps 
$(1,1)$, double-flatsteps $(2,0)$, and downsteps $(1,-1)$ that starts at the origin, never 
drops below the $x$-axis, and terminates on the $x$-axis. The 
terminal point necessarily has even $x$-coordinate, say $2n$; then we 
call it
a Schr\"{o}der $n$-path and we say its  size is $n$. 
The number of Schr\"{o}der $n$-paths is $r_{n}$. 
Here is a bijection from DF Dyck paths to  Schr\"{o}der paths 
that preserves size and sends marks to double-flatsteps.
Given a DF Dyck path, locate the upsteps starting at a marked 
vertex along with their matching downsteps. Then delete these 
upsteps (and their marks) and replace each matching downstep by a 
double-flatstep. An example with  3 marks and size 8 is illustrated. 
Note that a mark on the initial vertex of the Dyck path 
corresponds to the first double-flatstep at ground level in the 
Schr\"{o}der path; thus there are just as many Schr\"{o}der $n$-paths 
with a double-flatstep at ground level as without.
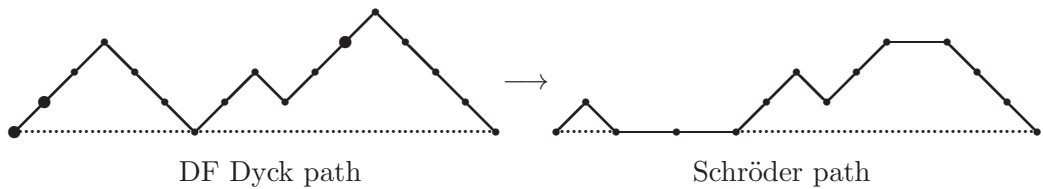
\Einheit=0.4cm
\[
\Pfad(-17,0),3334443343334444\endPfad
\SPfad(-17,0),1111111111111111\endSPfad
\Label\o{ \textrm{DF Dyck path}}(-8.5,-2)
\Label\o{ \textrm{Schr\"{o}der path}}(8.5,-2)
\Label\o{ \longrightarrow}(0,1)
\NormalPunkt(-17,0)
\NormalPunkt(-16,1)
\DuennPunkt(-15,2)
\DuennPunkt(-14,3)
\DuennPunkt(-13,2)
\DuennPunkt(-12,1)
\DuennPunkt(-11,0)
\DuennPunkt(-10,1)
\DuennPunkt(-9,2)
\DuennPunkt(-8,1)
\DuennPunkt(-7,2)
\NormalPunkt(-6,3)
\DuennPunkt(-5,4)
\DuennPunkt(-4,3)
\DuennPunkt(-3,2)
\DuennPunkt(-2,1)
\DuennPunkt(-1,0)
\DuennPunkt(1,0)
\DuennPunkt(2,1)
\DuennPunkt(3,0)
\DuennPunkt(5,0)
\DuennPunkt(7,0)
\DuennPunkt(8,1)
\DuennPunkt(9,2)
\DuennPunkt(10,1)
\DuennPunkt(11,2)
\DuennPunkt(12,3)
\DuennPunkt(14,3)
\DuennPunkt(15,2)
\DuennPunkt(16,1)
\DuennPunkt(17,0)
\Pfad(1,0),3411113343311444\endPfad
\SPfad(1,0),11\endSPfad
\SPfad(7,0),1111111111\endSPfad
\]

\newpage 
%\vspace*{7mm}

{\Large \textbf{4.\quad The DUtoDXD Bijection }  }

We define a bijection on Dyck $n$-paths with the following properties.
\begin{enumerate}
    \item  It sends \#\,$DU$s (valleys) to \#\,$DXD$s. Since 
    \#\,valleys has the Narayana distribution \cite{Narayana2}, so 
    does \#\,$DXD$s. Many other statistics 
    having the Narayana distribution on Dyck paths are given in \cite{Narayana2}.

    \item  It sends the paths with a terminal descent of even length 
    to the hill-free paths, thereby giving a bijection between two 
    manifestations of the Fine numbers \cite[p.\,263]{fine}. 

    \item  It sends the nonempty paths all of whose descents to ground level 
    have odd length to the paths that start $UD$, thereby giving a 
    bijective proof for item $j$ in \cite[Ex.\ 6.19]{ec2}.
\end{enumerate}

The bijection $\phi$ can be defined recursively as follows. First, 
$\phi(\eps)=\eps$ and we consider nonempty paths $P$ by the parity of the 
length of the terminal descent.

\noindent \ (i)\ terminal descent has even length. Here $P$ has the 
form (uniquely)
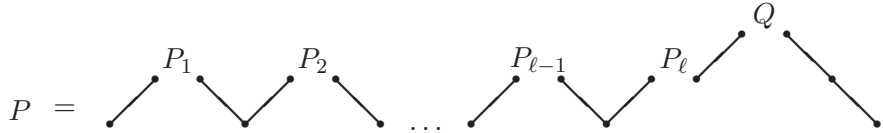
\Einheit=0.6cm
\[
\Pfad(-7,0),3\endPfad
\Pfad(-5,1),43\endPfad
\Pfad(-2,1),4\endPfad
\Pfad(1,0),3\endPfad
\Pfad(3,1),43\endPfad
\Pfad(6,1),3\endPfad
\Pfad(8,2),44\endPfad
%\SPfad(-6,0),1111111111111111\endSPfad
\DuennPunkt(-7,0)
\DuennPunkt(-6,1)
\DuennPunkt(-5,1)
\DuennPunkt(-4,0)
\DuennPunkt(-3,1)
\DuennPunkt(-2,1)
\DuennPunkt(-1,0)
\DuennPunkt(1,0)
\DuennPunkt(2,1)
\DuennPunkt(3,1)
\DuennPunkt(4,0)
\DuennPunkt(5,1)
\DuennPunkt(6,1)
\DuennPunkt(7,2)
\DuennPunkt(8,2)
\DuennPunkt(9,1)
\DuennPunkt(10,0)
\Label\o{P}(-9,-0.1)
\Label\o{=}(-8,-0.1)
\Label\o{P_{1}}(-5.5,1)
\Label\o{P_{2}}(-2.5,1)
\Label\o{\ldots}(0,-0.5)
\Label\o{P_{\ell-1}}(2.5,1)
\Label\o{P_{\ell}}(5.5,1)
\Label\o{Q}(7.5,2)
\]
with $\ell\ge 1,\ P_{1},\ldots,P_{\ell}$ arbitrary Dyck paths, and 
$Q$ a Dyck path (possibly empty) whose terminal descent has even 
length. Then 
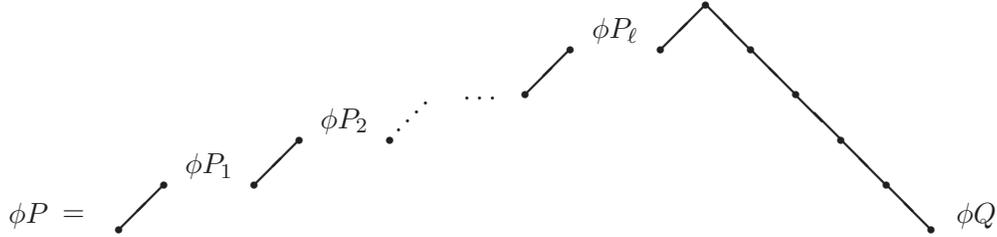
\Einheit=0.6cm
\[
\Pfad(-8,0),3\endPfad
\Pfad(-5,1),3\endPfad
\Pfad(1,3),3\endPfad
\Pfad(4,4),344444\endPfad
\SPfad(-2,2),3\endSPfad
\Label\o{\ldots}(0,2.5)
%\SPfad(-6,0),1111111111111111\endSPfad
\DuennPunkt(-8,0)
\DuennPunkt(-7,1)
\DuennPunkt(-5,1)
\DuennPunkt(-4,2)
\DuennPunkt(-2,2)
\DuennPunkt(1,3)
\DuennPunkt(2,4)
\DuennPunkt(4,4)
\DuennPunkt(5,5)
\DuennPunkt(6,4)
\DuennPunkt(7,3)
\DuennPunkt(8,2)
\DuennPunkt(9,1)
\DuennPunkt(10,0)
\Label\o{\phi P}(-10,-0.1)
\Label\o{=}(-9,-0.1)
\Label\o{\phi P_{1}}(-6,1)
\Label\o{\phi P_{2}}(-3,2)
\Label\o{\phi P_{\ell}}(3,4)
\Label\o{\phi Q}(11,-0.1)
\]

\noindent (ii)\ terminal descent has odd length. Here $P$ has the 
form
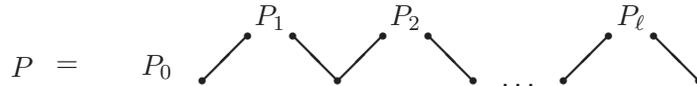
\Einheit=0.6cm
\[
\Pfad(-4,0),3\endPfad
\Pfad(-2,1),43\endPfad
\Pfad(1,1),4\endPfad
\Pfad(4,0),3\endPfad
\Pfad(6,1),4\endPfad
\Label\o{\ldots}(3,-0.5)
%\SPfad(-6,0),1111111111111111\endSPfad
\DuennPunkt(-4,0)
\DuennPunkt(-3,1)
\DuennPunkt(-2,1)
\DuennPunkt(-1,0)
\DuennPunkt(0,1)
\DuennPunkt(1,1)
\DuennPunkt(2,0)
\DuennPunkt(4,0)
\DuennPunkt(5,1)
\DuennPunkt(6,1)
\DuennPunkt(7,0)
\Label\o{ P}(-8,-0.1)
\Label\o{=}(-7,-0.1)
\Label\o{ P_{1}}(-2.5,1)
\Label\o{ P_{2}}(0.5,1)
\Label\o{ P_{\ell}}(5.5,1)
\Label\o{ P_{0}}(-5,-0.1)
\]
with  $\ell\ge 1,\ P_{0}, P_{1},\ldots,P_{\ell}$ Dyck paths (possibly empty) 
whose terminal descent has even length. Then 
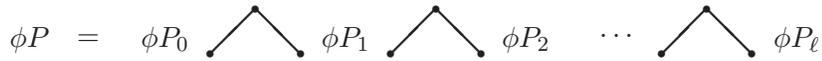
\Einheit=0.6cm
\[
\Pfad(-5,0),34\endPfad
\Pfad(-1,0),34\endPfad
\Pfad(5,0),34\endPfad
\Label\o{\ldots}(4,0)
%\SPfad(-6,0),1111111111111111\endSPfad
\DuennPunkt(-5,0)
\DuennPunkt(-4,1)
\DuennPunkt(-3,0)
\DuennPunkt(-1,0)
\DuennPunkt(0,1)
\DuennPunkt(1,0)
\DuennPunkt(5,0)
\DuennPunkt(6,1)
\DuennPunkt(7,0)
\Label\o{\phi P}(-9,-0.1)
\Label\o{=}(-7.7,-0.1)
\Label\o{\phi P_{0}}(-6,-0.1)
\Label\o{\phi P_{1}}(-2,-0.1)
\Label\o{\phi P_{2}}(2,-0.1)
\Label\o{\phi P_{\ell}}(8,-0.1)
\]
Note that all $UD$s terminating the path are preserved; in 
particular, $\phi$ preserves the property ``ends $DD$''. Clearly, 
$\phi$ sends Dyck paths whose terminal descent has odd length to 
paths containing at least one hill and, by induction, it sends paths
whose terminal descent has even length to hill-free paths. These 
facts allow the process to be reversed and yield a similar recursive 
definition of the inverse $\psi$. Again there are two cases:

\noindent \ (i) $P$ is nonempty hill-free. 
Here $P$ has the form
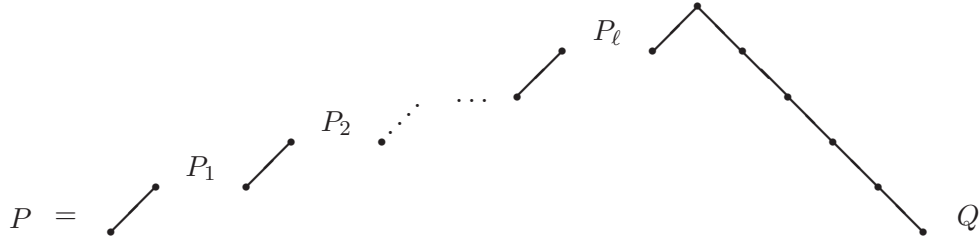
\Einheit=0.6cm
\[
\Pfad(-8,0),3\endPfad
\Pfad(-5,1),3\endPfad
\Pfad(1,3),3\endPfad
\Pfad(4,4),344444\endPfad
\SPfad(-2,2),3\endSPfad
\Label\o{\ldots}(0,2.5)
%\SPfad(-6,0),1111111111111111\endSPfad
\DuennPunkt(-8,0)
\DuennPunkt(-7,1)
\DuennPunkt(-5,1)
\DuennPunkt(-4,2)
\DuennPunkt(-2,2)
\DuennPunkt(1,3)
\DuennPunkt(2,4)
\DuennPunkt(4,4)
\DuennPunkt(5,5)
\DuennPunkt(6,4)
\DuennPunkt(7,3)
\DuennPunkt(8,2)
\DuennPunkt(9,1)
\DuennPunkt(10,0)
\Label\o{ P}(-10,-0.1)
\Label\o{=}(-9,-0.1)
\Label\o{ P_{1}}(-6,1)
\Label\o{ P_{2}}(-3,2)
\Label\o{ P_{\ell}}(3,4)
\Label\o{ Q}(11,-0.1)
\]
with $\ell \ge 1,\ P_{1},P_{2},\ldots,P_{\ell}$ arbitrary Dyck 
paths and $Q$ (possibly empty) a hill-free Dyck path. Then
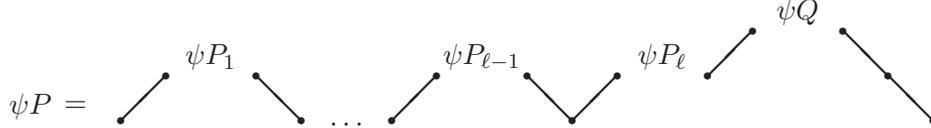
\Einheit=0.6cm
\[
\Pfad(-8,0),3\endPfad
\Pfad(-5,1),4\endPfad
\Pfad(-2,0),3\endPfad
\Pfad(1,1),43\endPfad
\Pfad(5,1),3\endPfad
\Pfad(8,2),44\endPfad
%\SPfad(-6,0),1111111111111111\endSPfad
\DuennPunkt(-8,0)
\DuennPunkt(-7,1)
\DuennPunkt(-5,1)
\DuennPunkt(-4,0)
\DuennPunkt(-2,0)
\DuennPunkt(-1,1)
\DuennPunkt(1,1)
\DuennPunkt(2,0)
\DuennPunkt(3,1)
\DuennPunkt(5,1)
\DuennPunkt(6,2)
\DuennPunkt(8,2)
\DuennPunkt(9,1)
\DuennPunkt(10,0)
\Label\o{\psi P}(-10,-0.1)
\Label\o{=}(-9,-0.1)
\Label\o{\psi P_{1}}(-6,1)
\Label\o{\ldots}(-3,-0.5)
\Label\o{\psi P_{\ell-1}}(0,1)
\Label\o{\psi P_{\ell}}(4,1)
\Label\o{\psi Q}(7,2)
\]
\noindent  (ii) $P$ contains hills. Here $P$ has the form

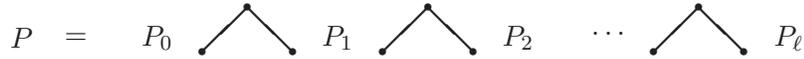
\Einheit=0.6cm
\[
\Pfad(-5,0),34\endPfad
\Pfad(-1,0),34\endPfad
\Pfad(5,0),34\endPfad
\Label\o{\ldots}(4,0)
%\SPfad(-6,0),1111111111111111\endSPfad
\DuennPunkt(-5,0)
\DuennPunkt(-4,1)
\DuennPunkt(-3,0)
\DuennPunkt(-1,0)
\DuennPunkt(0,1)
\DuennPunkt(1,0)
\DuennPunkt(5,0)
\DuennPunkt(6,1)
\DuennPunkt(7,0)
\Label\o{ P}(-9,-0.1)
\Label\o{=}(-7.8,-0.1)
\Label\o{ P_{0}}(-6,-0.1)
\Label\o{P_{1}}(-2,-0.1)
\Label\o{ P_{2}}(2,-0.1)
\Label\o{ P_{\ell}}(8,-0.1)
\]

with $\ell \ge 1,\ P_{0},P_{1},\ldots,P_{\ell}$ hill-free Dyck 
paths (possibly empty). Then

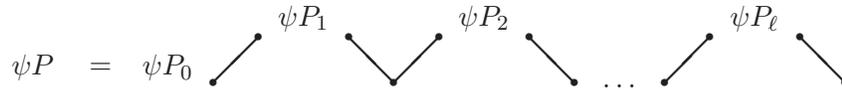
\Einheit=0.6cm
\[
\Pfad(-5,0),3\endPfad
\Pfad(-2,1),43\endPfad
\Pfad(2,1),4\endPfad
\Pfad(5,0),3\endPfad
\Pfad(8,1),4\endPfad
\Label\o{\ldots}(4,-0.5)
%\SPfad(-6,0),1111111111111111\endSPfad
\DuennPunkt(-5,0)
\DuennPunkt(-4,1)
\DuennPunkt(-2,1)
\DuennPunkt(-1,0)
\DuennPunkt(0,1)
\DuennPunkt(2,1)
\DuennPunkt(3,0)
\DuennPunkt(5,0)
\DuennPunkt(6,1)
\DuennPunkt(8,1)
\DuennPunkt(9,0)
\Label\o{\psi P}(-9,-0.1)
\Label\o{=}(-7.5,-0.1)
\Label\o{\psi  P_{0}}(-6,-0.1)
\Label\o{\psi  P_{1}}(-3,1)
\Label\o{\psi  P_{2}}(1,1)
\Label\o{\psi  P_{\ell}}(7,1)
\]

To show (by induction) that $\phi$ sends \#\,$DU$s to \#\,$DXD$s, 
let $\mu(P)$ denote the number of $DU$s in $P$ and $\nu(P)$ the number 
of $DXD$s. Then in case (i)---terminal descent has even length---we 
find
\[
\mu(P)=\ell-1 +\mu(P_{1})+\ldots+ \mu(P_{\ell-1}) +\mu(P_{\ell}) 
+\mu(Q) +[P_{\ell}\ne \eps].
\]
Here $[A]$ is the Iverson notation: $[A]=1$ if $A$ is true and $=0$ 
if $A$ is false. The term $[P_{\ell}\ne \eps]$ is present because the 
$U$ following $P_{\ell}$ is part of a $DU$ precisely when $P_{\ell}$ 
is nonempty. We also have
\[
\nu(\phi P)=\ell-1 +\nu(\phi P_{1})+\ldots+ \nu(\phi P_{\ell-1}) +\nu(\phi P_{\ell}) 
+\nu(\phi Q) +[\phi P_{\ell}\ne \eps].
\]
Here we have used the fact that $\phi Q$ does not start $UD$ because 
it is hill-free and the initial term $\ell-1$ records the $\ell-1 \ DDD$s 
contained in the first descent to ground level of $\phi P$. So 
induction yields $\mu(P)=\nu(\phi P)$ in case (i) and a similar 
argument works for case (ii). Because of this key property, we will 
subsequently refer to $\phi$ as the DUtoDXD bijection and to its 
inverse as DXDtoDU.

Next, we show that DXDtoDU sends the paths that start $UD$ to the 
paths all of whose descents to ground level have odd length. A path 
$P$ that starts $UD$ 
necessarily has $P_{0}=\eps$ in the form (ii) above for hill-containing paths 
and $\psi P$, as defined, has all descents to ground level of odd 
length because $\psi$ sends hill-free paths to paths whose terminal 
descent has even length. Similarly, $\psi$ sends all other paths to 
images that contain at least one even-length descent to ground level.

A recursively defined bijection is often easier to work with but it 
is also of some interest to have an explicit description. Here is a 
cut-and-paste description of DXDtoDU. Given a Dyck path, first color 
red the last $D$ of each $DDD$. Then color blue the middle $U$ of 
each $DUU$ unless its matching downstep is immediately followed by a 
red $D$.
\vspace*{2mm}
\Einheit=0.6cm
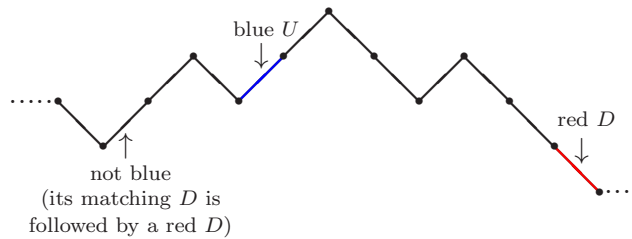
\[
\Pfad(-6,1),433433443444\endPfad
\blue{\Pfad(-2,1),3\endPfad}
\red{\Pfad(5,0),4\endPfad}
\SPfad(-7,1),1\endSPfad
\SPfad(6,-1),1\endSPfad
\DuennPunkt(-6,1)
\DuennPunkt(-5,0)
\DuennPunkt(-4,1)
\DuennPunkt(-3,2)
\DuennPunkt(-2,1)
\DuennPunkt(-1,2)
\DuennPunkt(0,3)
\DuennPunkt(1,2)
\DuennPunkt(2,1)
\DuennPunkt(3,2)
\DuennPunkt(6,-1)
\DuennPunkt(4,1)
\DuennPunkt(5,0)
\Label\o{\textrm{ {\scriptsize blue $U$}}}(-1.5,2.2)
\Label\o{\downarrow}(-1.5,1.6)
\Label\o{\textrm{ {\scriptsize red $D$}}}(5.6,0.2)
\Label\o{\downarrow}(5.6,-0.5)
\Label\o{\textrm{ {\scriptsize not blue}}}(-4.5,-1)
\Label\o{\textrm{ {\scriptsize (its matching $D$ is}}}(-4.5,-1.6)
\Label\o{\textrm{ {\scriptsize followed by a red $D$)}}}(-4.5,-2.2)
\Label\o{\uparrow}(-4.5,-0.4)
\]

\noindent Process the blue $U$s one at a time as follows (the resulting 
path is independent of the order of processing).
Delete the Dyck subpath starting at a blue $U$ and terminating at its 
matching $D$, and reinsert at the rightmost peak preceding it just 
before deletion. Then process the red $D$s
in the resulting path as follows (simultaneously, if you like). 
For each red $D$ locate the matching upstep $U$ 
of its predecessor step (the predecessor is necessarily a downstep); 
then delete each red $D$ and reinsert just before the corresponding 
$U$. An example is illustrated below with red steps labeled 
1,2,\ldots and blue steps labeled $a,b,\ldots$.

\noindent To invert the map we must first recover the red $D$s and blue $U$s in 
the image path. This is a little more involved. A \emph{valley} $D$ is 
one whose successor step is a $U$. The \emph{associated Dyck path} of 
an upstep in a Dyck path is the longest Dyck subpath starting at the 
upstep. The red $D$s are recovered as the valley $D$s for which the 
associated Dyck path of the successor upstep contains at least one 
descent of even length to (its own) ground level. Then the blue $U$s 
are recovered as the $U$s whose associated Dyck path (i) is immediately preceded 
by an upstep (i.e. $U$ is the second step of a double-rise), (ii) has 
terminal descent of even length, and (iii) is NOT immediately 
followed by a red $D$. Once the colored steps have been determined, 
the reverse cut-and-paste procedure is clear.

\Einheit=0.5cm
\[
\Pfad(-14,19),3333443344\endPfad
\SPfad(-14,19),1111111111111111111111111111\endSPfad
\red{\Pfad(-4,21),44\endPfad}
\red{\Pfad(9,20),4\endPfad}
\Label\o{ \textrm{{\scriptsize 1}}}(-3.4,20.3)
\Label\o{\textrm{{\scriptsize 2}}}(-2.4,19.3)
\Label\o{\textrm{{\scriptsize 3}}}(9.6,19.3)
\Label\o{ \textrm{{\scriptsize 1}}}(0.6,10.3)
\Label\o{\textrm{{\scriptsize 2}}}(1.6,9.3)
\Label\o{\textrm{{\scriptsize 3}}}(12.6,10.3)
\Label\o{ \textrm{{\scriptsize 1}}}(-6.4,0.3)
\Label\o{\textrm{{\scriptsize 2}}}(-12.4,0.3)
\Label\o{\textrm{{\scriptsize 3}}}(4.6,1.3)
\Label\o{ \textrm{{\scriptsize a}}}(-1.4,18.8)
\Label\o{\textrm{{\scriptsize b}}}(4.7,18.8)
\Label\o{\textrm{{\scriptsize c}}}(10.6,18.8)
\Label\o{ \textrm{{\scriptsize a}}}(-5.4,12.8)
\Label\o{\textrm{{\scriptsize b}}}(3.7,9.8)
\Label\o{\textrm{{\scriptsize c}}}(6.6,12.8)
\Label\o{ \textrm{{\scriptsize a}}}(-3.4,1.8)
\Label\o{\textrm{{\scriptsize b}}}(3.7,0.8)
\Label\o{\textrm{{\scriptsize c}}}(7.6,2.8)
\blue{\Pfad(-2,19),3\endPfad
\Pfad(4,19),3\endPfad
\Pfad(10,19),3\endPfad}
\Pfad(-1,20),34434\endPfad
\Pfad(5,20),3344\endPfad
\Pfad(11,20),344\endPfad
\Label\o{ \textrm{Dyck path $P$}}(0,24)
\Label\o{ \textrm{DXDtoDU$(P)$}}(0,5)
\Label\o{\textrm{{\scriptsize process\,}}\downarrow\textrm{{\scriptsize 
\,blues}}}(0,17)
%The default is	\PfadDicke=1pt.
\DuennPunkt(-14,19)
\DuennPunkt(-13,20)
\DuennPunkt(-12,21)
\DuennPunkt(-11,22)
\DuennPunkt(-10,23)
\DuennPunkt(-9,22)
\DuennPunkt(-8,21)
\DuennPunkt(-7,22)
\DuennPunkt(-6,23)
\DuennPunkt(-5,22)
\DuennPunkt(-4,21)
\DuennPunkt(-3,20)
\DuennPunkt(-2,19)
\DuennPunkt(-1,20)
\DuennPunkt(0,21)
\DuennPunkt(1,20)
\DuennPunkt(2,19)
\DuennPunkt(3,20)
\DuennPunkt(4,19)
\DuennPunkt(5,20)
\DuennPunkt(6,21)
\DuennPunkt(7,22)
\DuennPunkt(8,21)
\DuennPunkt(9,20)
\DuennPunkt(10,19)
\DuennPunkt(11,20)
\DuennPunkt(12,21)
\DuennPunkt(13,20)
\DuennPunkt(14,19)
\Pfad(-14,9),33334433\endPfad
\SPfad(-14,9),1111111111111111111111111111\endSPfad
\red{\Pfad(0,11),44\endPfad}
\red{\Pfad(12,11),4\endPfad}
\blue{\Pfad(-6,13),3\endPfad
\Pfad(3,10),3\endPfad
\Pfad(6,13),3\endPfad}
\Pfad(-5,14),34444\endPfad
\Pfad(2,9),3\endPfad
\Pfad(4,11),33\endPfad
\Pfad(7,14),34444\endPfad
\Pfad(13,10),4\endPfad
\Label\o{\textrm{{\scriptsize process\,}}\downarrow\textrm{{\scriptsize 
\,reds}}}(0,7)
%The default is	\PfadDicke=1pt.
\DuennPunkt(-14,9)
\DuennPunkt(-13,10)
\DuennPunkt(-12,11)
\DuennPunkt(-11,12)
\DuennPunkt(-10,13)
\DuennPunkt(-9,12)
\DuennPunkt(-8,11)
\DuennPunkt(-7,12)
\DuennPunkt(-6,13)
\DuennPunkt(-5,14)
\DuennPunkt(-4,15)
\DuennPunkt(-3,14)
\DuennPunkt(-2,13)
\DuennPunkt(-1,12)
\DuennPunkt(0,11)
\DuennPunkt(1,10)
\DuennPunkt(2,9)
\DuennPunkt(3,10)
\DuennPunkt(4,11)
\DuennPunkt(5,12)
\DuennPunkt(6,13)
\DuennPunkt(7,14)
\DuennPunkt(8,15)
\DuennPunkt(9,14)
\DuennPunkt(10,13)
\DuennPunkt(11,12)
\DuennPunkt(12,11)
\DuennPunkt(13,10)
\DuennPunkt(14,9)
\Pfad(-14,0),3\endPfad
\SPfad(-14,0),1111111111111111111111111111\endSPfad
\red{\Pfad(-13,1),4\endPfad
\Pfad(-7,1),4\endPfad}
\red{\Pfad(4,2),4\endPfad}
\blue{\Pfad(-4,2),3\endPfad
\Pfad(3,1),3\endPfad
\Pfad(7,3),3\endPfad}
\Pfad(-12,0),33344\endPfad
\Pfad(-6,0),33\endPfad
\Pfad(-3,3),344443\endPfad
\Pfad(5,1),33\endPfad
\Pfad(8,4),344444\endPfad
\DuennPunkt(-14,0)
\DuennPunkt(-13,1)
\DuennPunkt(-12,0)
\DuennPunkt(-11,1)
\DuennPunkt(-10,2)
\DuennPunkt(-9,3)
\DuennPunkt(-8,2)
\DuennPunkt(-7,1)
\DuennPunkt(-6,0)
\DuennPunkt(-5,1)
\DuennPunkt(-4,2)
\DuennPunkt(-3,3)
\DuennPunkt(-2,4)
\DuennPunkt(-1,3)
\DuennPunkt(0,2)
\DuennPunkt(1,1)
\DuennPunkt(2,0)
\DuennPunkt(3,1)
\DuennPunkt(4,2)
\DuennPunkt(5,1)
\DuennPunkt(6,2)
\DuennPunkt(7,3)
\DuennPunkt(8,4)
\DuennPunkt(9,5)
\DuennPunkt(10,4)
\DuennPunkt(11,3)
\DuennPunkt(12,2)
\DuennPunkt(13,1)
\DuennPunkt(14,0)
\]
\begin{center}
    An example of DXDtoDU
\end{center}

The properties of DUtoDXD that were proved above by induction from the 
recursive definition can also be seen directly from the explicit description.

\vspace*{7mm}

{\Large \textbf{5.\quad Review of Known Bijections }  }

For our main application of DUtoDXD, we need to review some well 
known bijections involving Dyck paths. The simplest bijection of all 
merely flips the path in the vertical, an involution that we call 
ReversePath. There is another involution $\phi$ on Dyck paths, due to 
Emeric Deutsch \cite{deutschInvol}, that can be defined recursively 
as follows.
\vspace*{2mm}
\Einheit=0.6cm
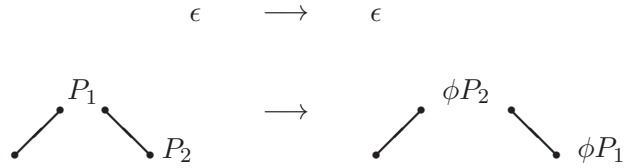
\[
\Pfad(-7,0),3\endPfad
\Pfad(-5,1),4\endPfad
\Pfad(1,0),3\endPfad
\Pfad(4,1),4\endPfad
\DuennPunkt(-7,0)
\DuennPunkt(-6,1)
\DuennPunkt(-5,1)
\DuennPunkt(-4,0)
\DuennPunkt(1,0)
\DuennPunkt(2,1)
\DuennPunkt(4,1)
\DuennPunkt(5,0)
\Label\o{ P_{1}}(-5.5,1)
\Label\o{\longrightarrow}(-1,0.5)
\Label\o{ P_{2}}(-3.4,-0.3)
\Label\o{\phi  P_{2}}(3,1)
\Label\o{ \phi P_{1}}(6,-0.3)
\Label\o{\longrightarrow}(-1,2.7)
\Label\o{\eps}(-3,2.7)
\Label\o{\eps}(1,2.7)
\]
The main properties of this involution  that we need 
are as follows (proof by induction): 

\ (i) it sends \#\,$UU$s to \#\,$DU$s (in fact, being an involution, it 
interchanges them)

(ii) it sends ``length of first descent'' to ``$1+\#\,UD$s that 
terminate the path''. In particular, it send paths with a short 
(resp. long) first descent to paths that end $DD$ (resp. $UD$).

We mention in passing that, using the well known ``walkaround'' bijection 
from binary trees to Dyck paths, $\phi$ has a simple explicit 
description: given a Dyck path, pass to the corresponding binary 
tree, flip it in the vertical, then pass back. We will refer to this 
involution as DeutschInvol.

A Levine-($r,s$) pair (apparently first considered in \cite{levine}\,) is a pair $(B,T)$ of 
nonintersecting lattice paths of steps North ($N=(0,1)$) and East 
($E=(1,0)$), where $B$ (the bottom  path) runs from (1,0) to 
$(r,s-1)$ and $T$ (the top path) runs from $(0,1)$ to $(r-1,s)$. 
There is a well known bijection, LevineToDyck, from Levine-$(r,s)$ pairs by 
way of ascent/descent sequences to Dyck-$(r+s-1)$ paths with $r$ 
peaks, illustrated  below with $r=4,s=5$.

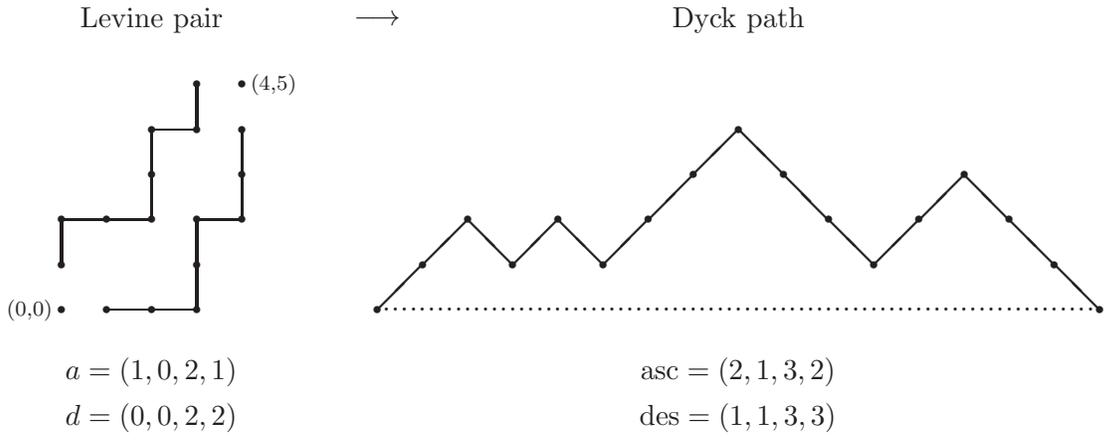
\Einheit=0.6cm
\[
\Pfad(-11,1),2112212\endPfad
\Pfad(-10,0),1122122\endPfad
\Pfad(-4,0),3343433344433444\endPfad
\SPfad(-4,0),1111111111111111\endSPfad
\Label\l{\textrm{{\scriptsize (0,0)}}}(-11.2,0)
\Label\r{\textrm{{\scriptsize (4,5)}}}(-6.8,5)
\DuennPunkt(-11,0)
\DuennPunkt(-7,5)
\DuennPunkt(-11,1)
\DuennPunkt(-11,2)
\DuennPunkt(-10,2)
\DuennPunkt(-9,2)
\DuennPunkt(-9,3)
\DuennPunkt(-9,4)
\DuennPunkt(-8,4)
\DuennPunkt(-8,5)
\DuennPunkt(-10,0)
\DuennPunkt(-9,0)
\DuennPunkt(-8,0)
\DuennPunkt(-8,1)
\DuennPunkt(-8,2)
\DuennPunkt(-7,2)
\DuennPunkt(-7,3)
\DuennPunkt(-7,4)
\DuennPunkt(-4,0)
\DuennPunkt(-3,1)
\DuennPunkt(-2,2)
\DuennPunkt(-1,1)
\DuennPunkt(0,2)
\DuennPunkt(1,1)
\DuennPunkt(2,2)
\DuennPunkt(3,3)
\DuennPunkt(4,4)
\DuennPunkt(5,3)
\DuennPunkt(6,2)
\DuennPunkt(7,1)
\DuennPunkt(8,2)
\DuennPunkt(9,3)
\DuennPunkt(10,2)
\DuennPunkt(11,1)
\DuennPunkt(12,0)
\Label\o{a=(1,0,2,1)}(-9,-1.8)
\Label\o{d=(0,0,2,2)}(-9,-2.8)
\Label\o{\textrm{asc} = (2,1,3,2)}(4,-1.8)
\Label\o{\textrm{des} = (1,1,3,3)}(4,-2.8)
\Label\o{\textrm{Levine pair}}(-9,6)
\Label\o{\textrm{Dyck path}}(4,6)
\Label\o{\longrightarrow}(-4,6)
\]

\vspace*{2mm}

\noindent Set $a=(a_{i})_{i=1}^{r}$ with $a_{i}$ = number of consecutive $N$s 
preceding the $i$th $E$ in the top path ($a_{r}=\#\,$terminal $N$s) 
and likewise define $d=(d_{i})_{i=1}^{r}$ for the bottom path. Then 
use $1+a$ and $1+d$ as the ascent lengths (left to right) and descent 
lengths respectively of the corresponding Dyck-$(r+s-1)$ path. The 
Dyck path has $r$ ascents and hence $r$ peaks. Note that if the top 
path $T$ ends $N$ (resp. $E$), then the last ascent of the Dyck path 
is long (resp. short).

\vspace*{7mm}

{\Large \textbf{6.\quad Dyck Paths by Long Interior Inclines}  }

The number of Dyck $n$-paths containing $k$ long \emph{nonterminal} inclines 
is $\frac{1}{n+1}\binom{n}{k}\binom{n}{k+1}$---the Narayana 
distribution \cite{Narayana2}.
The number of Dyck $n$-paths containing $k$ long \emph{interior} inclines is 
given in the following table for small $n,k$,
\[
\begin{array}{c|cccccc}
	n^{\textstyle{\,\backslash \,k}}&0 &1 &2 & 3&4 & 5\\
\hline 
	1&    1 &   & & &  & \\
 	2&    2 &  & & & &  \\
	3&    3 & 2  & & & & \\
	4&    4 & 8  &  2 & & &  \\ 
	5&    5 & 20 & 15 &2 & &  \\
	6&    6 & 40 & 60 &24 &2 &  \\
	7&    7 & 70 &175 &140 &35 &2  \\

 \end{array}
\]
and we wish to show that the $(n,k)$ entry is
$\frac{2}{n+1}\binom{n+1}{k+2} \binom{n-2}{k}$, thus establishing a 
combinatorial interpretation of the identity
\[
C_{n}=\frac{2}{n+1}\sum_{k=0}^{n-2}\binom{n+1}{k+2} 
\binom{n-2}{k}.
\]

We have $\frac{2}{n+1}\binom{n+1}{k+2} 
\binom{n-2}{k}=\binom{n-2}{k}\binom{n}{k+1}-\binom{n-2}{k-1}\binom{n}{k+2}$ 
and, by the beautiful Gessel-Viennot determinant theorem for nonintersecting 
paths \cite{Gessel-Viennot,thebook3}, 
this is the number of pairs $(B,T)$ of nonintersecting lattice paths 
consisting of steps $N$ and $E$ where $B$ runs from $(1,0)$ to 
$(k+2,n-1-k)$ and $T$ runs from $(0,1)$ to $(k,n-1-k)$. On this set, 
the map ``delete last step of $B$, flip it ($N\leftrightarrow E$), then 
append it to $T$'' is a bijection (i) from the pairs in which $B$ 
ends $N$ to the set $\a$ of Levine-$(k+2,n-1-k)$ pairs in which the 
top path ends $E$, and (ii) from the pairs in which $B$ 
ends $E$ to the set $\b$ of Levine-$(k+1,n-k)$ pairs in which the 
top path ends $N$.

We are about to exhibit a sequence of bijections from \a (resp. \b) 
to Dyck $n$-paths containing $k$ long interior inclines and ending 
$DD$ (resp. $UD$). But first we need to relate the number of long interior 
inclines in a Dyck path $P$ to the number of $DXD$s in the elevated 
path $E(P)$ (to elevate a path, prepend $U$ and append $D$, thus 
$E(P)=UPD$).

\textbf{Lemma}\quad
   \emph{ For $n\ge 1$ and $P$ a Dyck $n$-path,
  \#\, long interior inclines in $P +\ \#\,DXD$s in $E(P)=n-1$}.
  
\textbf{Proof}\quad In $E(P)$, each of its $n+1\ D$s is either 
(i) preceded by $UU$, (ii) preceded by $UD$, (iii) an interior $D$ in 
a long descent, or (iv) the last $D$ in a long descent. The $D$s in 
the respective classes 
correspond in an obvious way to (i) long ascents, (ii) $DUD$s, 
(iii) $DDD$s, and (iv) long descents in $E(P)$. Since the first and 
last inclines in $E(P)$ are necessarily long, the result follows. \qed 

Since each noninitial $U$ in a Dyck $n$-path $P$ is preceded either by 
a $U$ or a $D$, we also have the obvious result for $P$: \#\,$UU$s  
$+\,\#\,DU$s $=n-1$. 
And it is all but obvious that \#\,$UD$s $=1+\#\,DU$s.

Finally, here are the promised sequences of bijections:
\begin{eqnarray*}
    \a & \xrightarrow{\text{LevineToDyck}} & \textrm{Dyck $n$-paths 
    with a short last ascent and $k+2$ peaks}  \\
     &  \xrightarrow{\text{ReversePath}} &  \textrm{Dyck $n$-paths 
    with a short first descent and $n-2-k\ UU$s} \\
     & \xrightarrow{\text{DeutschInvol}} & \textrm{Dyck $n$-paths that 
     end $DD$ with $n-2-k\ DU$s}  \\
     & \xrightarrow{\text{DUtoDXD}} & \textrm{Dyck $n$-paths that 
     end $DD$ with $n-2-k\ DXD$s}  \\
     &  & \textrm{= Dyck $n$-paths that end $DD$ with $n-1-k\ DXD$s in 
     $E(P)$}  \\
     &  & \textrm{= Dyck $n$-paths with $k$ long interior 
     inclines that end $DD$. }
\end{eqnarray*}

Similarly,
\begin{eqnarray*}
    \b & \xrightarrow{\text{LevineToDyck}} & \textrm{Dyck $n$-paths 
    with a long last ascent and $k+1$ peaks}  \\
     &  \xrightarrow{\text{ReversePath}} &  \textrm{Dyck $n$-paths 
    with a long first descent and $n-1-k\ UU$s} \\
     & \xrightarrow{\text{DeutschInvol}} & \textrm{Dyck $n$-paths that 
     end $UD$ with $n-1-k\ DU$s}  \\
     & \xrightarrow{\text{DUtoDXD}} & \textrm{Dyck $n$-paths that 
     end $UD$ with $n-1-k\ DXD$s}  \\
     &  & \textrm{= Dyck $n$-paths that end $UD$ with $n-1-k\ DXD$s in 
     $E(P)$}  \\
     &  & \textrm{= Dyck $n$-paths with $k$ long interior 
     inclines that end $UD$. }
\end{eqnarray*}

\vspace*{7mm}

{\Large \textbf{7.\quad Another Catalan Number Identity }  }

An identity related to that of the preceding section is
\[
\hspace*{30mm}C_{n}=\sum_{k=0}^{n-1}\left(\binom{n-1}{k}^{2}-
\binom{n+1}{k+2}\binom{n-3}{k-2}\right) \hspace*{35mm}n\ge 1
\]
Again using Gessel-Viennot, the summand is the number of pairs $(B,T)$ 
of nonintersecting lattice paths, $B$ running from $(2,0)$ to 
$(k+2,n-1-k)$ and $T$ running from $(0,0)$ to $(k,n-1-k)$.
Delete the last step of $B$, flip it and append to $T$ as in the last 
section. Likewise, delete the first step of $T$, flip it and prepend to 
$B$. This produces Levine pairs and then LevineToDyck is a bijection 
to Dyck $n$-paths, establishing the identity. By methods similar to 
those of the last section, the identity in fact counts Dyck-$n$ paths 
by the statistic $X:=X_{1}+X_{2}$ where $X_{1}= $ \#\,long interior inclines 
and $X_{2}=$ [path is strict] where strict means the path returns to 
ground level only at its terminal point. Is there a more ``natural'' 
statistic on Dyck paths (or on some other manifestation of the Catalan 
numbers) with this distribution?

\vspace*{7mm}

{\Large \textbf{8 \quad Ordered Trees by Nodes Adjacent to a Leaf}  }

Following Deutsch and Shapiro \cite{2Motzkin}), we define a node in an ordered tree 
to be a vertex that is neither a leaf nor the root. Thus the
vertices are partitioned into 3 classes: root, nodes, leaves.

The number of ordered trees on $n$ edges containing $k$ nodes adjacent 
to a leaf is given in the following table for small $n,k$,
\[
\begin{array}{c|ccccc}
	n^{\textstyle{\,\backslash \,k}} & 0 & 1 & 2 & 3 & 4  \\
\hline 
	1&    1 &   & & & \\
 	2&    1 & 1& & &   \\
	3&    1 & 4  & & &  \\
	4&    1 & 10  & 3 & &   \\ 
	5&    1 & 20 & 21 & &   \\
	6&    1 & 35 & 84 & 12 &   \\
	7&    1 & 56 &252 & 120&   \\
        8&    1 & 84 &630 & 660& 55   \\

 \end{array}
\]
and we will show that the $(n,k)$ entry is $\frac{1}{n+1}\sum_{k}\binom{n+1}{2k+1} \binom{n+k}{k}$.  
First, we show directly that $\binom{n+1}{2k+1} \binom{n+k}{k}$  
counts certain marked balanced paths. Then we exhibit a bijection from these 
marked paths to (plain) balanced $n$-paths with $n-2k$ odd-length ascents. 
Next we apply \S 3 to see that 
$\frac{1}{n+1}\sum_{k}\binom{n+1}{2k+1} \binom{n+k}{k}$ counts Dyck 
$n$-paths with $n-2k$ odd-length ascents. Finally, we exhibit a 
bijection from the latter paths to Dyck paths that correspond, under 
a familiar bijection (variously known as the walkaround, glove or accordion 
bijection), to
ordered trees on $n$ edges containing $k$ nodes adjacent to a leaf.
Incidentally, on full binary trees the statistic ``\# nodes adjacent to 
a leaf'' has the Touchard distribution \cite{superCatalan}.

Recall DF  refers to paths with downstep-free vertices available for 
marking. Now we show that $\binom{n+1}{2k+1} 
\binom{n+k}{k}$ counts DF balanced $n$-paths containing 
$n-2k$ marked vertices and satisfying the property: for each ascent, its length 
minus its number of marked vertices is even.
Start with a row of $n$ upsteps. They determine $n+1$ gaps: the $n-1$ 
gaps between them and a gap at each end. Choose $2k+1$ of these 
gaps---$\binom{n+1}{2k+1}$ choices---and label them 
$G_{1},G_{2},\ldots,G_{2k+1}$ from left to right. Insert a mark in 
each of the remaining $n-2k$ gaps. Distribute $n$ downsteps 
arbitrarily among every other labeled gap, that is, among gaps 
$G_{1},G_{3},G_{5},\ldots,G_{2k+1}$---$\binom{n+k}{k}$ choices---and 
concatenate to form a typical path of the type specified. An example 
with $n=8,\,k=2$ is illustrated (slanted lines represent the upsteps).

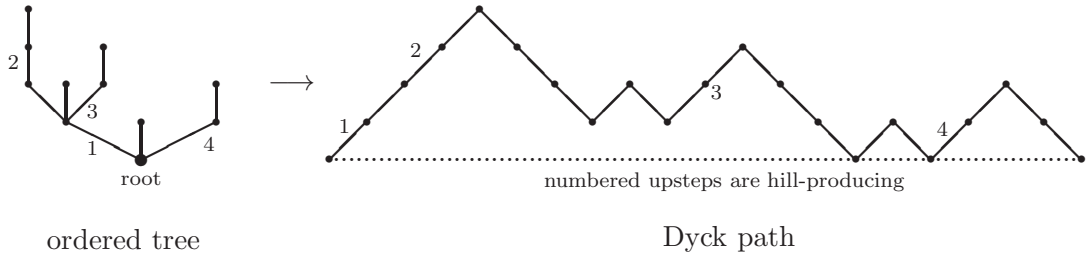
\Einheit=0.5cm
\[
\Pfad(-8,0),3\endPfad
\Pfad(-6,0),3\endPfad
\Pfad(-4,0),3\endPfad
\Pfad(-2,0),3\endPfad
\Pfad(0,0),3\endPfad
\Pfad(2,0),3\endPfad
\Pfad(4,0),3\endPfad
\Pfad(6,0),3\endPfad
\Pfad(8,0),3\endPfad
\Label\o{G_{1}}(-8.5,0)
\Label\o{G_{2}}(-6.5,0)
\Label\o{G_{3}}(-4.5,0)
\Label\o{G_{4}}(3.5,0)
\Label\o{G_{5}}(5.5,0)
\NormalPunkt(-2.5,0.5)
\NormalPunkt(-0.5,0.5)
\NormalPunkt(1.5,0.5)
\NormalPunkt(7.5,0.5)
\]

\noindent Say $3\ D$s go into gap $G_{1},\ 1 \ D$ goes into $G_{3}$ and 
$4\ D$s go into 
$G_{5}$. (In general, of course, some of these gaps might remain empty.) Then the resulting DFV-marked path is 
\[
\Pfad(-8,0),4443343333344443\endPfad
\SPfad(-8,0),1111111111111111\endSPfad
\NormalPunkt(-1,-1)
\NormalPunkt(0,0)
\NormalPunkt(1,1)
\NormalPunkt(8,0)
\DuennPunkt(-8,0)
\DuennPunkt(-7,-1)
\DuennPunkt(-6,-2)
\DuennPunkt(-5,-3)
\DuennPunkt(-4,-2)
\DuennPunkt(-3,-1)
\DuennPunkt(-2,-2)
\DuennPunkt(2,2)
\DuennPunkt(3,3)
\DuennPunkt(4,2)
\DuennPunkt(5,1)
\DuennPunkt(6,0)
\DuennPunkt(7,-1)
\]

\noindent Here is a bijection from these marked paths to (plain) 
balanced $n$-paths, that sends marks to odd-length ascents. Given such 
a marked path, first erase the marks (if any) on vertices at ground 
level. Each vertex $V$ above ground level has a matching vertex in 
the path: the first one encountered directly east of $V$. For each 
marked vertex $V$ above ground level, delete the upstep starting at 
$V$ and reinsert it at $V$'s matching vertex (without the mark). Do likewise for marked 
vertices below ground level, replacing the word ``east'' by ``west''. 
The path illustrated above yields (transferred steps shown in green merely as a 
visual aid)
\[
\Pfad(-8,3),44433\endPfad
\Pfad(-2,3),433344\endPfad
\Pfad(5,4),443\endPfad
\SPfad(-8,3),1111111111111111\endSPfad
\green{\Pfad(-3,2),3\endPfad
\Pfad(4,3),3\endPfad}
\DuennPunkt(-8,3)
\DuennPunkt(-7,2)
\DuennPunkt(-6,1)
\DuennPunkt(-5,0)
\DuennPunkt(-4,1)
\DuennPunkt(-3,2)
\DuennPunkt(-2,3)
\DuennPunkt(-1,2)
\DuennPunkt(0,3)
\DuennPunkt(1,4)
\DuennPunkt(2,5)
\DuennPunkt(3,4)
\DuennPunkt(4,3)
\DuennPunkt(5,4)
\DuennPunkt(6,3)
\DuennPunkt(7,2)
\DuennPunkt(8,3)
\]

\noindent The transferred steps above ground level appear in the image 
path as the initial steps of odd-length ascents  above ground level 
that do not start the path. Similarly, transferred steps below ground level appear in the image 
path as the terminal steps of odd-length ascents  below ground level 
that do not end the path. The erased marks at ground level correspond 
to odd-length ascents that either cross ground level or 
start or end the path. Hence this map sends ``\#\,marks'' to 
``\#\,odd-length ascents'' and is invertible.
From \S3 we 
conclude that $\frac{1}{n+1}\binom{n+1}{2k+1} \binom{n+k}{k}$ counts 
Dyck $n$-paths with $n-2k$ odd-length ascents.

An ordered tree on $n$ edges corresponds to a Dyck $n$-path under an 
obvious bijection: burrow up all edges from the root and open out 
accordion-style so that each edge produces an upstep and matching 
downstep, as illustrated.

\Einheit=0.5cm
\[
\Pfad(-14,2),22\endPfad
\Pfad(-14,2),46\endPfad
\Pfad(-11,0),52\endPfad
\Pfad(-13,1),2\endPfad
\Pfad(-13,1),32\endPfad
\Pfad(-11,0),2\endPfad
\Pfad(-6,0),33334443433444343344\endPfad
\SPfad(-6,0),11111111111111111111\endSPfad
\Label\o{\longrightarrow}(-7,1.5)
\Label\u{ \textrm{{\scriptsize root}}}(-11,0.1)
\Label\u{ \textrm{{\scriptsize numbered upsteps are 
hill-producing}}}(4.5,0)
\NormalPunkt(-11,0)
\DuennPunkt(-14,4)
\DuennPunkt(-14,3)
\DuennPunkt(-14,2)
\DuennPunkt(-11,1)
\DuennPunkt(-12,2)
\DuennPunkt(-12,3)
\DuennPunkt(-13,1)
\DuennPunkt(-13,2)
\DuennPunkt(-9,1)
\DuennPunkt(-9,2)
\DuennPunkt(-6,0)
\DuennPunkt(-5,1)
\DuennPunkt(-4,2)
\DuennPunkt(-3,3)
\DuennPunkt(-2,4)
\DuennPunkt(-1,3)
\DuennPunkt(-0,2)
\DuennPunkt(1,1)
\DuennPunkt(2,2)
\DuennPunkt(3,1)
\DuennPunkt(4,2)
\DuennPunkt(5,3)
\DuennPunkt(6,2)
\DuennPunkt(7,1)
\DuennPunkt(8,0)
\DuennPunkt(9,1)
\DuennPunkt(10,0)
\DuennPunkt(11,1)
\DuennPunkt(12,2)
\DuennPunkt(13,1)
\DuennPunkt(14,0)
\Label\u{ \textrm{ ordered tree}}(-11.6,-1.5)
\Label\u{ \textrm{ Dyck path}}(4.5,-1.5)
\Label\o{ \textrm{{\scriptsize 1}}}(-12.3,-0.2)
\Label\o{ \textrm{{\scriptsize 2}}}(-14.4,2.1)
\Label\o{ \textrm{{\scriptsize 3}}}(-12.3,0.8)
\Label\o{ \textrm{{\scriptsize 4}}}(-9.2,-0.1)
\Label\o{ \textrm{{\scriptsize 1}}}(-5.6,.4)
\Label\o{ \textrm{{\scriptsize 2}}}(-3.7,2.4)
\Label\o{ \textrm{{\scriptsize 3}}}(4.3,1.3)
\Label\o{ \textrm{{\scriptsize 4}}}(10.3,0.3)
\]

\noindent Each non-root vertex in the tree has a unique parent 
edge: the first one on the path to the root. The parent edges of 
nodes (\,= non-root non-leaf vertices) that are adjacent to a leaf, numbered in the illustration, correspond to 
upsteps $U$ in the Dyck path with the following property: the 
Dyck subpath lying strictly between $U$ and its matching $D$ contains 
at least one hill. We'll call such a $U$ a \emph{hill-producing} upstep.

Now consider placing nonoverlapping dimers (dominos), each the length 
of two steps, to cover as many upsteps as possible in a Dyck path. 
They can all be covered except for one upstep in each odd-length ascent. 
Clearly, the number of odd-length ascents in a Dyck $n$-path has the 
same parity (even or odd) as $n$. Hence, if this number is $n-2k$, 
then $k$ dimers can be placed. Thus we have shown that the number 
of Dyck $n$-paths accommodating $k$ dimers on the upsteps is 
$\frac{1}{n+1}\binom{n+1}{2k+1} \binom{n+k}{k}$. In the next section we exhibit 
a bijection  on Dyck paths that sends
\#\,accommodated dimers  to \#\,hill-producing upsteps. As we have 
seen, the latter correspond to nodes adjacent to a leaf in an ordered 
tree, and the result follows.

\vspace*{7mm}

{\Large \textbf{9 \quad The Dimer to Hill-Producing Upstep Bijection    }  }

Here we define a bijection, $\phi$, on Dyck $n$-paths that sends
\#\,accommodated dimers  to \#\,hill-producing upsteps.
Our definition of $\phi$ is recursive and it is convenient to 
introduce a little more notation. A strict Dyck path is a 
Dyck path with exactly one return to ground level (necessarily at the 
end of the path). In particular, a strict Dyck path is nonempty. The 
interior of a strict Dyck path is the Dyck path obtained by deleting 
its initial upstep and terminal downstep. Every nonempty Dyck path is 
representable uniquely as a concatenation of one or more strict Dyck 
paths, called its components. We abbreviate first ascent length as
FAL. The interior of the first component 
of a nonempty Dyck path $P$ is denoted $I(P)$, and $I^{2}(P)=I(I(P))$ 
and so on. Thus if $P$ has FAL $=k$, then $I^{j}(P)\ne \eps$ for $1 
\le j \le k-1$ and $I^{k}(P)=\eps$. Recall a Dyck path is Fine if it 
contains no hills (peaks at level 1); otherwise it is Hill. 
%otherwise we say it is a Hill path.

Now define $\phi\, \eps =\eps,\ \phi\, UD=UD$, and for nonstrict $P$ with 
components $P_{1},\ldots,P_{k}\ (k\ge 2),\ \phi P =\phi P_{1}\,\phi 
P_{2}\,\ldots\,\phi P_{k}$. It remains to define $\phi$ on strict Dyck 
paths with FAL $\ge 2$ and the image is another strict Dyck 
path of the same size. In the Table below, the left column partitions 
strict Dyck paths with FAL $\ge 2$ into six classes (four with FAL 
even and two with FAL odd). The right column lists properties that (as 
we will see) distinguish the image classes. The diagrams following the 
Table define $\phi$ on each of the six classes in terms of $\phi$ 
acting on smaller paths. 

\newcommand{\rb}[1]{\raisebox{-1.5ex}[0pt]{#1}}

\vspace*{5mm}

%\begin{table}
\centerline{
\begin{tabular}{|c||c|c|} \hline
 \rb{Case}  & \multicolumn{2}{|c|}{\rule[-4mm]{0mm}{10mm}$\phi\ :\ P \longrightarrow\ Q$} \\ \hline\hline
                         & \underline{FAL of $P$ is even and:\vphantom{Q}}    & {\rule[-4mm]{0mm}{10mm}\underline{$I(Q)$ is Hill and:}}  \\ 
{\rule[-3mm]{0mm}{8mm}1.} & $I(P)$ is Fine                                 & $I^{2}(Q)$ is Hill       \\ 
{\rule[-3mm]{0mm}{8mm}2.} & $I(P)$ is Hill,\ $I^{2}(P)=\eps$               & $I^{2}(Q)=\eps$          \\ 
{\rule[-3mm]{0mm}{8mm}3.} & $I(P)$ is Hill,\ $I^{2}(P)\ne \eps$ is Fine    & $I^{2}(Q)\ne\eps$ is Fine, $I^{3}(Q)$ is Hill    \\ 
{\rule[-3mm]{0mm}{8mm}4.} & $I(P)$ is Hill,\ $I^{2}(P)\ne \eps$ is Hill    & $I^{2}(Q)\ne\eps$ is Fine, $I^{3}(Q)$ is Fine       \\ \hline \hline
                         & \underline{FAL of $P$ is odd  and:\vphantom{Q}}     & {\rule[-4mm]{0mm}{10mm}\underline{$I(Q)$ is nonempty Fine and:}}  \\ 
{\rule[-3mm]{0mm}{8mm}5.} & $I(P)$ is Fine                                 & $I^{2}(Q)$ is Hill                              \\ 
{\rule[-3mm]{0mm}{8mm}6.} & $I(P)$ is Hill                                 & $I^{2}(Q)$ is Fine                              \\ \hline
\end{tabular}}
%\caption{Classification of strict Dyck paths for the bijection $\phi$}
%\end{table}
\vspace*{-3mm}
\begin{center}
    Classification of strict Dyck paths for the bijection $\phi$ \\
{\footnotesize  [FAL = first ascent length,  $I(P)$ denotes interior of first 
component of $P$] }   
%    $I(P)$ is interior of first component of $P$
\end{center}

\vspace*{2mm}

\noindent We note in passing that, restricted to the paths in the top 
two boxes, $\phi$ is a bijection between two known manifestations of 
the Fine numbers \cite[p.\,263]{fine}.

In the diagrammed paths, $P_{1},P_{2},\ldots$ represent Dyck paths subject only to the restrictions specified.

\noindent Case 1.
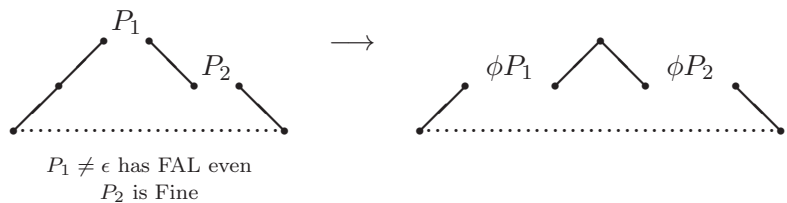
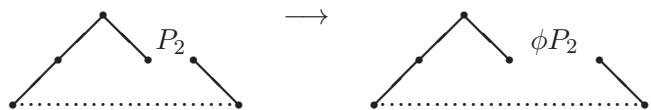
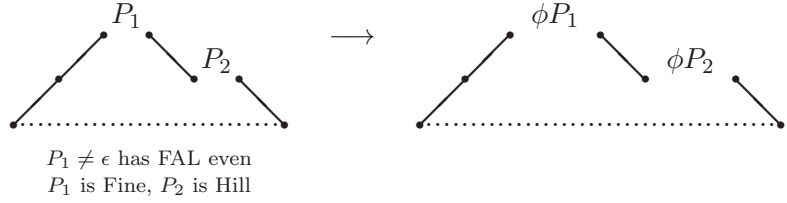
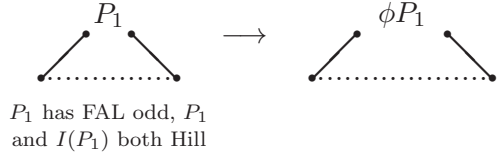
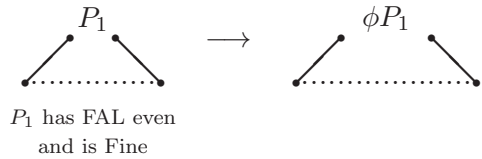
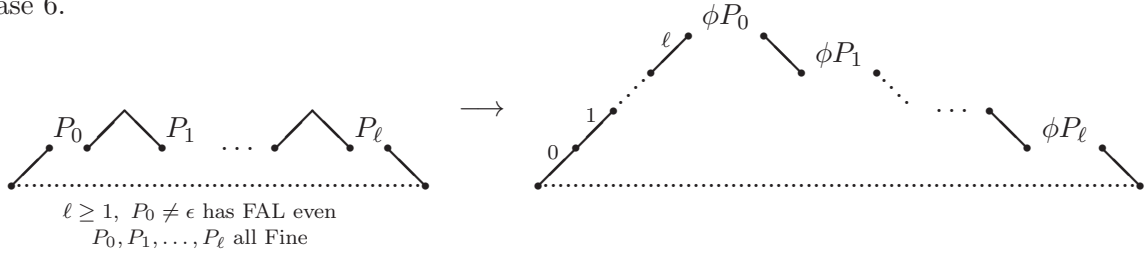
\Einheit=0.6cm
\[
\Pfad(-9,0),33\endPfad
\Pfad(-6,2),4\endPfad
\Pfad(-4,1),4\endPfad
\Pfad(0,0),3\endPfad
\Pfad(3,1),34\endPfad
\Pfad(7,1),4\endPfad
\SPfad(-9,0),111111\endSPfad
\SPfad(0,0),11111111\endSPfad
\Label\o{P_{1}}(-6.5,2)
\Label\o{P_{2}}(-4.5,1)
\Label\o{\longrightarrow}(-1.5,1.5)
\Label\o{\phi P_{1}}(2,1)
\Label\o{\phi P_{2}}(6,1)
\DuennPunkt(-9,0)
\DuennPunkt(-8,1)
\DuennPunkt(-7,2)
\DuennPunkt(-6,2)
\DuennPunkt(-5,1)
\DuennPunkt(-4,1)
\DuennPunkt(-3,0)
\DuennPunkt(0,0)
\DuennPunkt(1,1)
\DuennPunkt(3,1)
\DuennPunkt(4,2)
\DuennPunkt(5,1)
\DuennPunkt(7,1)
\DuennPunkt(8,0)
\Label\o{\textrm{{\scriptsize $P_{1}\ne \eps$  has FAL even}}}(-6,-1.2)
\Label\o{\textrm{{\scriptsize $P_{2}$  is Fine}}}(-6,-1.8)
\]

\noindent Case 2.
\[
\Pfad(-7,0),334\endPfad
\Pfad(-3,1),4\endPfad
\Pfad(1,0),334\endPfad
\Pfad(6,1),4\endPfad
\SPfad(-7,0),11111\endSPfad
\SPfad(1,0),111111\endSPfad
\Label\o{P_{2}}(-3.5,1)
\Label\o{\longrightarrow}(-0.5,1.5)
\Label\o{\phi P_{2}}(5,1)
\DuennPunkt(-7,0)
\DuennPunkt(-6,1)
\DuennPunkt(-5,2)
\DuennPunkt(-4,1)
\DuennPunkt(-3,1)
\DuennPunkt(-2,0)
\DuennPunkt(1,0)
\DuennPunkt(2,1)
\DuennPunkt(3,2)
\DuennPunkt(4,1)
\DuennPunkt(6,1)
\DuennPunkt(7,0)
\]

\noindent Case 3.
\[
\Pfad(-9,0),33\endPfad
\Pfad(-6,2),4\endPfad
\Pfad(-4,1),4\endPfad
\Pfad(0,0),33\endPfad
\Pfad(4,2),4\endPfad
\Pfad(7,1),4\endPfad
\SPfad(-9,0),111111\endSPfad
\SPfad(0,0),11111111\endSPfad
\Label\o{P_{1}}(-6.5,2)
\Label\o{P_{2}}(-4.5,1)
\Label\o{\longrightarrow}(-1.5,1.5)
\Label\o{\phi P_{1}}(3,2)
\Label\o{\phi P_{2}}(6,1)
\DuennPunkt(-9,0)
\DuennPunkt(-8,1)
\DuennPunkt(-7,2)
\DuennPunkt(-6,2)
\DuennPunkt(-5,1)
\DuennPunkt(-4,1)
\DuennPunkt(-3,0)
\DuennPunkt(0,0)
\DuennPunkt(1,1)
\DuennPunkt(2,2)
\DuennPunkt(4,2)
\DuennPunkt(5,1)
\DuennPunkt(7,1)
\DuennPunkt(8,0)
\Label\o{\textrm{{\scriptsize $P_{1}\ne \eps$  has FAL even}}}(-6,-1.2)
\Label\o{\textrm{{\scriptsize $P_{1}$  is Fine, $P_{2}$  is Hill}}}(-6,-1.8)
\]

\noindent Case 4.
\[
\Pfad(-5,0),3\endPfad
\Pfad(-3,1),4\endPfad
\Pfad(1,0),3\endPfad
\Pfad(4,1),4\endPfad
\SPfad(-5,0),111\endSPfad
\SPfad(1,0),1111\endSPfad
\Label\o{P_{1}}(-3.5,1)
\Label\o{\longrightarrow}(-0.5,0.5)
\Label\o{\phi P_{1}}(3,1)
\DuennPunkt(-5,0)
\DuennPunkt(-4,1)
\DuennPunkt(-3,1)
\DuennPunkt(-2,0)
\DuennPunkt(1,0)
\DuennPunkt(2,1)
\DuennPunkt(4,1)
\DuennPunkt(5,0)
\Label\o{\textrm{{\scriptsize $P_{1}$  has FAL odd, $P_{1}$}}}(-3.5,-1.2)
\Label\o{\textrm{{\scriptsize   and $I(P_{1})$  both Hill}}}(-3.5,-1.8)
\]

\noindent Case 5.
\[
\Pfad(-5,0),3\endPfad
\Pfad(-3,1),4\endPfad
\Pfad(1,0),3\endPfad
\Pfad(4,1),4\endPfad
\SPfad(-5,0),111\endSPfad
\SPfad(1,0),1111\endSPfad
\Label\o{P_{1}}(-3.5,1)
\Label\o{\longrightarrow}(-0.5,0.5)
\Label\o{\phi P_{1}}(3,1)
\DuennPunkt(-5,0)
\DuennPunkt(-4,1)
\DuennPunkt(-3,1)
\DuennPunkt(-2,0)
\DuennPunkt(1,0)
\DuennPunkt(2,1)
\DuennPunkt(4,1)
\DuennPunkt(5,0)
\Label\o{\textrm{{\scriptsize $P_{1}$  has FAL even}}}(-3.5,-1.2)
\Label\o{\textrm{{\scriptsize   and is Fine}}}(-3.5,-1.8)
\]

\noindent Case 6.
\Einheit=0.5cm
\[
\Pfad(-15,0),3\endPfad
\Pfad(-13,1),34\endPfad
\Pfad(-8,1),34\endPfad
\Pfad(-5,1),4\endPfad
\Pfad(-1,0),33\endPfad
\SPfad(1,2),3\endSPfad
\Pfad(2,3),3\endPfad
\SPfad(8,3),4\endSPfad
\Pfad(5,4),4\endPfad
\Pfad(11,2),4\endPfad
\Pfad(14,1),4\endPfad
\SPfad(-15,0),11111111111\endSPfad
\SPfad(-1,0),1111111111111111\endSPfad
\Label\o{\ldots}(-9,0.5)
\Label\o{\ldots}(10,1.5)
\Label\o{P_{0}}(-13.5,1)
\Label\o{P_{1}}(-10.5,1)
\Label\o{P_{\ell}}(-5.5,1)
\Label\o{\longrightarrow}(-2.5,1.5)
\Label\o{\phi P_{0}}(4,4)
\Label\o{\phi P_{1}}(7,3)
\Label\o{\phi P_{\ell}}(13,1)
\Label\o{\textrm{{\scriptsize $0$}}}(-0.6,0.4)
\Label\o{\textrm{{\scriptsize $1$}}}(0.4,1.4)
\Label\o{\textrm{{\scriptsize $\ell$}}}(2.4,3.4)
\DuennPunkt(-15,0)
\DuennPunkt(-14,1)
\DuennPunkt(-13,1)
\DuennPunkt(-11,1)
\DuennPunkt(-8,1)
\DuennPunkt(-6,1)
\DuennPunkt(-5,1)
\DuennPunkt(-4,0)
\DuennPunkt(-1,0)
\DuennPunkt(0,1)
\DuennPunkt(1,2)
\DuennPunkt(2,3)
\DuennPunkt(3,4)
\DuennPunkt(5,4)
\DuennPunkt(6,3)
\DuennPunkt(8,3)
\DuennPunkt(11,2)
\DuennPunkt(12,1)
\DuennPunkt(14,1)
\DuennPunkt(15,0)
\Label\o{\textrm{{\scriptsize $\ell\ge 1,\ P_{0}\ne \eps$ has FAL 
even}}}(-10,-1.2)
\Label\o{\textrm{{\scriptsize $P_{0},P_{1},\ldots,P_{\ell}$ all 
Fine}}}(-10,-1.9)
\]

By induction, $\phi$ sends paths with FAL even (resp. odd) to paths 
whose first component has an interior with (resp. without) hills, and 
the classification of image paths in the Table follows. Only in Case 
1 is the form of the image path not obviously unique and the hill displayed in that 
case is recovered as the \emph{last} hill in the interior of the 
image since $\phi$ sends Fine paths to Fine paths.
So an inverse for $\phi$ can be recursively defined by reversing 
everything.

To establish the key property of $\phi$, let $\nu(P) =\#\,$dimers accommodated 
by $P$, and $\mu(P)=\#\,$ hill-producing upsteps in $P$. It is easy to 
verify by induction that $\nu(P)=\mu(\phi P)$. For example, in Case 
1, $\nu(P) =1\textrm{ (due to initial $UU$)}+\nu(P_{1})+\nu(P_{2})$ while 
$\mu(\phi P) = 1\textrm{ (the first $U$ is hill-producing)} + \mu(\phi 
P_{1}) + \mu(\phi P_{2})$ and induction establishes this case. The 
others are similar.
This completes the proof.

%\vspace{7mm}
\newpage

{\Large \textbf{10 \quad Touchard for the Fine Numbers}  }

Touchard's Catalan number identity, $C_{n}=\sum_{k\ge 
0}\binom{n-1}{2k}2^{n-1-2k}C_{k}$, counts Dyck $n$-paths by number $k$ of 
long noninitial ascents. Here we show that the identity 
$F_{n}=\frac{1}{n+1}\sum_{k\ge 0}\binom{n-2-k}{k}$ $2^{n-2-2k}
\binom{n+1}{k+1}$ counts Fine $n$-paths by 
the same statistic.

The initial ascent of a nonempty Fine path is always long; so it 
suffices to show that $\frac{1}{n+1}\binom{n-1-k}{k-1}\binom{n-2k}{j}\binom{n+1}{k}$ is 
the number of Fine $n$-paths with $j$ short ascents and
$k$ long ascents, and then sum over $j$. To show this, with IA again short 
for interior ascent, define an 
\emph{IA balanced $(n,j,k)$-path} to be a balanced $n$-path 
with $k$ long ascents and $j$ marked IA vertices. An 
IA balanced path is \emph{Fine-like} if it contains no short 
ascents (and thus every ascent has at least one interior vertex) and 
the first interior vertex of each ascent is \emph{not} marked.

The number of Fine-like IA balanced $(n,j,k)$-paths is 
$\binom{n-1-k}{k-1}\binom{n-2k}{j}\binom{n+1}{k}$ as follows: start 
with a row of $n$ upsteps. From the $n-1$ gaps between them, choose 
$k$ nonconsecutive gaps including the first one---$\binom{n-k-1}{k-1}$ 
choices---to serve as the first interior vertices of the ascents. 
Choose $j$ gaps for marks from the $n-2k$ gaps available (the $k$ 
originally chosen gaps and their immediate predecessors are not 
available)---$\binom{n-2k}{j}$ choices. Finally, distribute $n\ D$s 
among the $k$ gaps preceding the originally chosen gaps and the 
``gap'' after the last upstep, with at least one $D$  in each gap 
save possibly the first and last of them---$\binom{n+1}{k}$ 
choices---to get a typical path of the type being counted.

Again from \S3 we have that there are $\frac{1}{n+1}\binom{n-1-k}{k-1}\binom{n-2k}{j}\binom{n+1}{k}$ 
Fine-like IA \emph{Dyck} $(n,j,k)$-paths. Now we exhibit a 
bijection from Fine-like IA Dyck paths to (plain) Fine paths that 
preserves size of path and number of long ascents, and turns marks 
into short ascents. An example with $n=12$ upsteps, $j=3$ marks, and 
$k=4$ long ascents is illustrated below.
For each marked vertex $v$, take the upstep $U$ terminating at $v$, 
locate the matching downstep $D$ for $U$, and transfer $U$ to the 
initial vertex of $D$ (erasing the mark) as indicated by the arrows 
in the top figure.
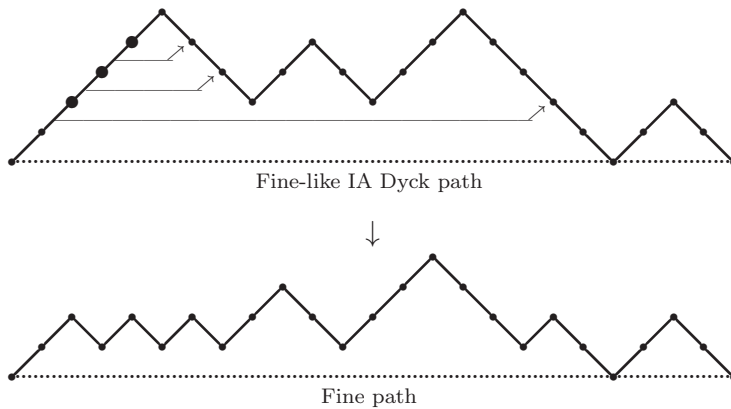
\Einheit=0.4cm
\[
\Pfad(-12,0),333334443344333444443344\endPfad
\SPfad(-12,0),111111111111111111111111\endSPfad
%\Label\o{\longrightarrow}(13,3)
\Label\o{\downarrow}(0,-3)
\Label\o{\textrm{{\scriptsize Fine-like IA Dyck path }}}(0,-1.3)
\DuennPunkt(-12,0)
\DuennPunkt(-11,1)
\NormalPunkt(-10,2)
\Label\o{\textrm{------}}(-7.6,2.7)
\Label\o{{\textrm{ \tiny $\nearrow$ }} }(-6.6,3.0)
\NormalPunkt(-9,3)
\Label\o{\textrm{------------}}(-7.6,1.7)
\Label\o{{\textrm{ \tiny $\nearrow$ }} }(-5.6,2.0)
\Label\o{\textrm{--------------------------------------------------}}(-2.65,0.7)
\Label\o{{\textrm{ \tiny $\nearrow$ }} }(5.4,1.0)
\NormalPunkt(-8,4)
\DuennPunkt(-7,5)
\DuennPunkt(-6,4)
\DuennPunkt(-5,3)
\DuennPunkt(-4,2)
\DuennPunkt(-3,3)
\DuennPunkt(-2,4)
\DuennPunkt(-1,3)
\DuennPunkt(0,2)
\DuennPunkt(1,3)
\DuennPunkt(2,4)
\DuennPunkt(3,5)
\DuennPunkt(4,4)
\DuennPunkt(5,3)
\DuennPunkt(6,2)
\DuennPunkt(7,1)
\DuennPunkt(8,0)
\DuennPunkt(9,1)
\DuennPunkt(10,2)
\DuennPunkt(11,1)
\DuennPunkt(12,0)
\]
\[
\Pfad(-12,0),334343433443334443443344\endPfad
\SPfad(-12,0),111111111111111111111111\endSPfad
\Label\o{\textrm{{\scriptsize Fine path }}}(0,-1.3)
\DuennPunkt(-12,0)
\DuennPunkt(-11,1)
\DuennPunkt(-10,2)
\DuennPunkt(-9,1)
\DuennPunkt(-8,2)
\DuennPunkt(-7,1)
\DuennPunkt(-6,2)
\DuennPunkt(-5,1)
\DuennPunkt(-4,2)
\DuennPunkt(-3,3)
\DuennPunkt(-2,2)
\DuennPunkt(-1,1)
\DuennPunkt(0,2)
\DuennPunkt(1,3)
\DuennPunkt(2,4)
\DuennPunkt(3,3)
\DuennPunkt(4,2)
\DuennPunkt(5,1)
\DuennPunkt(6,2)
\DuennPunkt(7,1)
\DuennPunkt(8,0)
\DuennPunkt(9,1)
\DuennPunkt(10,2)
\DuennPunkt(11,1)
\DuennPunkt(12,0)
\]

The transferred steps can be recovered as the short ascents
in the image path. To reverse the map, delete each short ascent $U$ 
and reinsert it with a mark on its top vertex at the initial vertex 
of the associated upstep (as defined in \S 2) of the downstep preceding 
$U$. Note that, as in the 
example, several upsteps may be reinserted at the same vertex. This 
completes the proof.

\end{document}